\documentclass[intlimits,twoside,a4paper]{article}

\usepackage{amsmath,amssymb}
\usepackage{graphicx}
\usepackage{wrapfig}
\usepackage[T2A]{fontenc}
\usepackage[cp1251]{inputenc}
\usepackage[eqsecnum]{cmpj}

\newtheorem{theorem}{Theorem}

\newtheorem{corollary}{Corollary}

\newtheorem{lemma}{Lemma}

\newtheorem{rmk}{Remark}

\issue{2011}{14}{2}{23003}

\doinumber{10.5488/CMP.14.23003}

\title[Phase transitions for the Potts model]%
{On phase transitions of the Potts model with three competing interactions on Cayley tree}

\author[H. Akin, S. Temir]{H.~Akin\refaddr{label1}\thanks{E-mail:
hasanakin69@gmail.com}\,,
S.~Temir\refaddr{label2}\thanks{E-mail: temirseyit@harran.edu.tr}}

\authorcopyright{H. Akin, S. Temir, 2011}

\addresses{
\addr{label1} Department of Mathematics, Faculty of Education,
Zirve University, 27260~Gaziantep, Turkey \addr{label2} Department
of Mathematics, Arts and Science Faculty Harran University,
63200~Sanliurfa, Turkey}

\date{Received November 15, 2010, in final form March 22, 2011}%

\begin{document}

\maketitle

\begin{abstract}
In the present paper we study a phase transition problem for the
Potts model with three competing interactions, the nearest
neighbors, the second neighbors and triples of neighbors and
non-zero external field on Cayley tree of order two. We prove that
for some parameter values of the model there is phase transition.
We reduce the problem of describing by limiting Gibbs measures to
the problem of solving a system of nonlinear functional equations.
We extend the results obtained by Ganikhodjaev and Rozikov [Math.
Phys. Anal. Geom., 2009, \textbf{12}, No.~2, 141--156] on phase
transition for the Ising model to the Potts model setting.
\pacs 05.50.+q, 64.60.-i, 64.60.De, 75.10.Hk
\end{abstract}
\section{Introduction}
It is well known that lattice spin system is a large class of
systems considered in statistical mechanic. It is also well known
that the structure of the lattice plays an important role in
examining the spin systems. The Potts model, which was introduced
by Potts (in 1952) as a generalization of the Ising model,  has
drawn an increased attention to its theoretical and experimental
aspects in recent years~\cite{W}. The idea originated from the
representation of the Ising model as interacting spins which can
be either parallel or antiparallel.

We consider a semi-infinite Cayley tree J$^{k}$ of order
$k\geqslant 2$, i.e., a graph without cycles with $(k+1)$ edges
issuing from each vertex except $x^{0}$ and with $k$ edges issuing
from the vertex $x^{0}$, which is called the tree root.  We
generalize the problem in~\cite{GAT} for the Potts model with
three competing interactions and non-zero external field. This
problem has been investigated for the Potts model with only two
competing interactions~\cite{GAT, GTA1}.

There are many approaches to the derivation of an equation or a
system of equations describing the limiting Gibbs measure for the
lattice models on Cayley tree (see~\cite{W,AT,B,P,Sinai,SC,S} for
 details). One of these approaches is based on recursive equations
for partition functions~\cite{W}. Another approach is Markov
random fields (see~\cite{GPW,K}). Naturally, both approaches lead
to the same equation~\cite{GPW}.

In this paper we employ the recursive equations for partition
functions. We study the problem of phase transition for the Potts
model with three competing interactions, the nearest neighbors,
the second neighbors and triples of neighbors and non-zero
external field on Cayley tree of order two. In~\cite{GR}, it has
been proved that there is a phase transition for some parameter
values of the Ising model. We generalize this problem for the
Potts model. In order to construct the Hamiltonian equation, we
consider the generalized Kroneker's $\delta$ function which well
agrees with theory of quadratic stochastic operators.

The paper is organized as follows. In section 2 we
present basic definitions and notations. In section 3 we present the
recursive equations for partition functions. Section 4 contains
the solutions of the system of nonlinear equations. These solutions
describe the presence or the absence of a phase transition.
Finally, the conclusions are presented in section~5.

\section{Preliminaries}
\textbf{Cayley tree.} The Cayley tree $\Gamma^{k}$ (see~\cite{B})
of order $k\geqslant 1$ is an infinite tree, i.e., a graph without
cycles, from each vertex of which exactly $k+1$ edges issue. Let
$\Gamma^{k}=(V, L,i)$, where $V$ is the set of vertices of
$\Gamma^{k}$, $L$ is the set of edges of $\Gamma^{k}$ and $i$ is
the incidence function associating each edge $\ell\in L$ with its
end points $x, y \in V.$ If $i(\ell)=\{x, y\}$, then $x$ and $y$
are called the nearest neighboring vertices and we write
$\ell=\langle x, y\rangle$. For $x,\  y \in V$, the distance $d(x,
y)$ on Cayley tree is defined by the formula
\begin{eqnarray*}
d(x, y)&=&\min \{d|x=x_{0}\,, x_{1}\,, x_{2}\,, \ldots, x_{d-1}\,,
x_{d}=y \in V
\ \text{such  that  the pairs} \\
&& \langle x_0\,, x_1\rangle, \ldots, \langle x_{d-1}\,, x_d
\rangle \text{ are neighboring vertices}\}.
\end{eqnarray*}
For the fixed $x^{0}\in V$ we have
\begin{eqnarray*}
W_{n}&=&\left\{x\in V: d(x^{0},x)=n\right\},
\\V_{n}&=&\left\{x\in V: d(x^{0},x)\leqslant  n\right\},
\\L_{n}&=&\left\{\ell= \langle x,y \rangle\in L:x,y\in V_{n} \right\}.
\end{eqnarray*}
A collection of the pairs  $\langle x_0\,, x_1 \rangle, \ldots,
\langle x_{d-1}\,, y \rangle$ is called a path from $x$ to $y$. We
write $x<y$ if the path from $x^{0}$ to $y$ goes through $x$. We
call the vertex $y$ a direct successor of $x$, if $y>x$ and $x,\
y$ are the nearest neighbors. The set of direct successors of $x$
is denoted by $S(x),$ i.e.,
\[
S(x)=\left\{y\in W_{n+1}\,:d(x,y)=1\right\}, \qquad x\in W_{n}\,.
\]
It is clear that for any vertex $x\neq x^{0}$, $x$ has $k$ direct
successors and $ x^{0}$ has $k+1$. The vertices $x$ and $y$ are
called the second neighbor which is denoted by $\rangle x, y
\langle$\,, if there exists a vertex $z\in V$  such that $x, z$
and $y, z$ are the nearest neighbors.

The second neighbors are called one-level neighbors, if vertices $x$ and $y$ belong to $W_{n}$ for some n, that
is, if they are located on the same level. We will consider only one-level
second neighbors. Three vertices $x, y$
and $z$ are called a triple of neighbors, which is denoted by
$\left\langle {x, y, z} \right\rangle $, if $\left\langle {x,y}
\right\rangle $ and $\left\langle {y, z} \right\rangle $ are nearest
neighbors and $x\neq z$.

For neighbors triple $\left\langle {x, y, z}
\right\rangle$ we consider the following Kroneker's symbol
\begin{equation}\label{eq0}
\delta _{\sigma (x)\sigma (y)\sigma (z)} =\left\{
{{\begin{array}{*{20}c}
 {1,} & \text{if}\ \  {\sigma (x)=\sigma (y)=\sigma (z)},&\\
 {{1}/{2},}& \text{if}\ \ {\sigma (x)=\sigma (y)\ne \sigma (z)}& \text{or} & {\sigma (x)\ne \sigma (y)=\sigma (z)},\\
 {0,}& \text{otherwise},
\end{array}}} \right.
\end{equation}
where $x, z\in W_{n}$ and $y\in W_{n-1}$ for some $n$. We can
write \eqref{eq0} as follows:
\[
\delta _{\sigma (x)\sigma (y)\sigma (z)} = \frac{1}{2}\left(\delta
_{\sigma (x)\sigma (y)}+\delta _{\sigma (y)\sigma (z)} \right).
\]

The generalized Kroneker's function in~\eqref{eq0} which
corresponds to the quadratic stochastic operator is an identity
transformation~\cite{MTA}.

\textbf{Model.} In the Potts model, the spin variables $\sigma(x)$
take their values on a discrete set $\Phi=\{1, 2,\dots , q\},\
q>2$ which are associated with each vertex of the tree $J^{k}$
(see~\cite{GTA} for details).

The Potts model with three competing interactions and external
field is defined by the following Hamiltonian
\begin{eqnarray}\label{hm}
H(\sigma )=-J\sum\limits_{\left\langle {x,y}
\right\rangle } {\delta _{\sigma (x)\sigma
(y)}-J_{1}\sum\limits_{\rangle x,y \langle}  } {\delta _{\sigma (x)\sigma
(y)}}-J_{2}\sum\limits_{\left\langle {x,y,z} \right\rangle } {\delta
_{\sigma (x)\sigma (y)\sigma (z)}}-h\sum\limits_{x\in V} {\delta
_{0\sigma (x)}}\,,
\end{eqnarray}
where the first sum ranges within all nearest neighbors, the
second sum ranges within all the second neighbors, the third sum
ranges within all triples of neighbors and the spin variables
$\sigma(x)$ may take the values 1, 2, 3. $J$, $J_{1}$\,, $J_{2}$
$\in\mathbb{R}$ are coupling constants and $h$ is external field.

In~\cite{GTA}, the authors have already obtained the exact
solutions of the Potts model with competing ternary and binary
interactions and external field on Cayley tree described by means
of the Hamiltonian~\eqref{hm}, where $h=0$ and $J_{1}=0$. In this
paper, we assume
that each parameter
$J,\ J_{1},\ J_{2},\ h$ is different from zero.

\section{Recursive equations for Partition Functions}
Let $\Lambda $ be a finite subset of V. We will denote the
restriction of $\sigma $ to $\Lambda$  by $\sigma (\Lambda)$. Let
$\bar \sigma $(V$\backslash \Lambda )$ be a fixed boundary
configuration. The total energy of $\sigma (\Lambda )$ under
condition $\bar \sigma $(V$\backslash \Lambda )$ is defined as
\begin{eqnarray}\label{hm2}
\nonumber H_{\Lambda}({\sigma(\Lambda)|{\bar \sigma } (V\backslash
\Lambda )} )=&-&J\sum\limits_{\langle x,y \rangle: x,y\in \Lambda
}{\delta _{\sigma (x)\sigma (y)}}-J\sum\limits_{\langle x,y
\rangle:{x\in \Lambda ,y\notin \Lambda }}\delta _{\sigma (x)\sigma
(y)}
\nonumber\\
&-&J_{1}\sum\limits_{\rangle x,y \langle:x,y\in \Lambda}\delta _{\sigma
(x)\sigma (y)}-J_{1}\sum\limits_{\rangle x,y \langle:x\in \Lambda ,y\notin
\Lambda}\delta _{\sigma (x)\sigma (y)}\nonumber\\
&-&J_{2}\sum\limits_{\langle x,y,z \rangle:x,y,z\in \Lambda}\delta _{\sigma
(x)\sigma (y)\sigma (z)}-J_{2}\sum\limits_{\langle x,y,z \rangle:x,z\notin
\Lambda,y\in \Lambda}\delta _{\sigma (x)\sigma (y)\sigma (z)}\nonumber\\
&-&h\sum\limits_{x\in \Lambda }\delta _{0\sigma (x)}.
\end{eqnarray}
Then, partition function $Z_{\Lambda} \left(\bar \sigma
(V\setminus \Lambda)\right)$ in volume $\Lambda $ under boundary
condition $\bar \sigma $(V$\backslash \Lambda )$ is defined as
\begin{equation}\label{par1}
Z_{\Lambda}({\bar \sigma (V\backslash \Lambda )})= \sum
\limits_{\sigma (\Lambda )\in \Omega (\Lambda )} \exp \left[
-\beta H_{\Lambda}(\sigma (\Lambda )| \bar \sigma (V\backslash
\Lambda ))\right],
\end{equation}
where $\Omega (\Lambda )$ is the set of all configuration in
volume $\Lambda $, and $\beta ={1}/{kT}$ is the inverse
temperature.

Instead of the configurations $\sigma $(V$_{n})$ and the partition
functions $Z_{V_n }$ in volume V$_{n}$ we use notations  $\sigma
_{n}$ and Z$^{(n)}$, respectively~\cite{GPW,GTA}.

It is well known that the Gibbs measure is a probability measure
used in many problems of probability theory and of statistical
mechanics~\cite{Ge}. The Gibbs measure of a configuration
$\sigma_\Lambda$ is defined by
\begin{equation*}
 \mu(\sigma_\Lambda)=\frac{\exp \left[-\beta H(\sigma_\Lambda)\right]}{Z_\Lambda(\sigma)}\,,
\end{equation*}
where $Z_\Lambda(\sigma)$ is the partition function.

As usual, one can introduce the notion of the Gibbs measure
(phase) of  the Potts model with a competing interaction on the
Cayley tree (see~\cite{P, Sinai, Ge} for details).

Let us decompose the
partition function Z$^{(n)}$ into the following summands
\[
Z^{(n)}=\sum\limits_{i=1}^q {Z_i^{(n)} },
\]
where
\begin{equation}\label{par2}
 Z_{i}^{(n)} =\sum\limits_{\sigma _n \in \Omega (V_{n} ):\sigma
(x^{0})=i}\exp \left[ -\beta H_{V_{n}} (\sigma _{n} | \bar \sigma
(V \backslash V_{n}))\right].
\end{equation}

Now the partition function is iteratively calculated by moving
from the boundary towards to the interior of the Cayley tree of
order two. The partial partition function of a branch on $n$
generations, where the innermost site is in state $i$, is denoted
by $Z_{i}^{(n)}$.

Let $S(x^0)=\{x^{1},\ x^{2} \}$, $\sigma(x^{0})=i$,
$\sigma(x^{1})=j$ and $\sigma(x^{2})=m$. For a given $n$, let us
define the partial partition functions of the next generation by
\[
Z_{i}^{(n)}=\overset{3}{\underset{j,m=1}{\sum }}\exp \left[\beta
J(\delta _{ij}+\delta _{im})+\beta J_{1}\delta _{jm}+\beta
J_{2}\delta _{jim}+h\delta _{1i}\right]Z_{j}^{(n-1)}Z_{m}^{(n-1)},
\]
where $i=1,\ 2,\ 3$.

Let $\theta = \exp(\beta J)$, $\theta _{1}= \exp(\beta J_{1})$,
$\theta _{2}= \exp(\beta J_{2})$ and $\theta_{3} = \exp(\beta h)$.
From~\eqref{hm2} and~\eqref{par1}, we can calculate the partial
partition functions as follows:
\begin{eqnarray*}
Z_1^{(n)}&=&\theta _3 \left\{\theta^{2} \theta _1\theta
_2\left(Z_1^{(n - 1)} \right)^2 + 2\theta \sqrt
{\theta_{2}}\left(Z_1^{(n - 1)} Z_2^{(n -
1)}+Z_1^{(n - 1)} Z_3^{(n - 1)}\right)\right.\\
 &+& \left.\theta _1\left[\left(Z_2^{(n - 1)} \right)^2
 +\left(Z_3^{(n - 1)} \right)^2\right]+2Z_2^{(n - 1)} Z_3^{(n -
 1)}\right\},\\
Z_2^{(n)}&=&\theta _1 \left(Z_1^{(n - 1)} \right)^2  + 2\theta\sqrt
{\theta_{2}} Z_2^{(n - 1)}\left(Z_1^{(n - 1)}+Z_3^{(n - 1)}\right)\\
&+& 2Z_1^{(n - 1)} Z_3^{(n - 1)}  + \theta ^2 \theta _1\theta
_2\left(Z_2^{(n - 1)}\right)^2+\theta _1 \left(Z_3^{(n - 1)} \right)^2,\\
Z_3^{(n)}&=& \theta _1\left(Z_1^{(n - 1)} \right)^2  + 2\theta\sqrt
{\theta_{2}} Z_3^{(n - 1)}\left(Z_1^{(n - 1)}+Z_2^{(n - 1)}\right)\\
&+&2Z_1^{(n - 1)} Z_2^{(n -
1)}+\theta^2\theta_1\theta_2\left(Z_3^{(n - 1)}\right)^2+ \theta
_1\left(Z_2^{(n - 1)} \right)^2.
\end{eqnarray*}
After replacement $u_n = {{Z_2^{(n)} }}/{{Z_1^{(n)} }}$ and $v_n =
{{Z_3^{(n)} }}/{{Z_1^{(n)} }}$ we have the following system of
recurrent equations
\begin{eqnarray}\label{dizi}
\theta _3 u_n &=& \frac{{\theta _1 + 2\theta\sqrt {\theta_{2}} u_{n
- 1}(1+v_{n - 1})  + 2v_{n - 1} + \theta^2\theta_1\theta_2u_{n -
1}^2+\theta _1 v_{n - 1}^2
}}{{\theta^2\theta_1\theta_2+2\theta\sqrt\theta _2(u_{n - 1} +v_{n
- 1})+\theta _1 (u_{n - 1}^2  +v_{n - 1}^2)+ 2u_{n - 1} v_{n -
1}}}\,,\nonumber\\[0.5ex]
\theta _3 v_n & =& \frac{{\theta _1  + 2u_{n - 1}  + 2\theta \sqrt
{\theta _2 } v_{n - 1} (1 + u_{n - 1} ) + \theta_1 u_{n - 1}^2  +
\theta ^2 \theta _1 \theta _2 v_{n - 1}^2
}}{\theta^2\theta_1\theta_2+2\theta\sqrt\theta _2(u_{n - 1} +v_{n
- 1})+\theta _1 (u_{n - 1}^2  +v_{n - 1}^2)+ 2u_{n - 1} v_{n - 1}}\,.
\end{eqnarray}
We define  the transformation
\begin{eqnarray}\label{fixedfunc}
 \textbf{F}=(F_1\,,F_2):
\mathbb{R}^2\rightarrow \mathbb{R}^2
\end{eqnarray}
with $u_n=F_1(u_{n-1}\,,v_{n-1})$ and
$v_n=F_2(u_{n-1}\,,v_{n-1})$. The fixed points of the
function~\eqref{fixedfunc} $\textbf{w}=\textbf{F}(\textbf{w})$,
where $\textbf{w}=(u_n\,,v_n)$, describe translation-invariant
Gibbs measures (phases) of the Potts model defined by the
Hamiltonian~\eqref{hm}.

The recursive equations~\eqref{dizi} can be written as
$w_n=\textbf{F}(w_{n-1})$, $n>0$. In the theory of dynamical
systems, $w_n$ is called trajectory of the initial point $w_1$
under the action of the mapping $\textbf{F}$. So, the asymptotic
behavior of $Z^{(n)}$ for $n\rightarrow \infty$ can be determined
by value of $\lim_nw_n$\,.

If $u=\lim_nu_{n}$ and $v=\lim_n v_{n}$ then we have the following equations;
\begin{eqnarray}\label{eq4}
\nonumber&& \theta _3 u  = \frac{{\theta _1 + 2\theta\sqrt
{\theta_{2}} u(1+v) + 2v + \theta^2\theta_1\theta_2u^2+\theta _1
v^2 }}{{\theta^2\theta_1\theta_2+2\theta\sqrt\theta _2(u
+v)+\theta _1
(u^2+v^2)+ 2uv}}\,,\\[0.5ex]
&&\theta _3 v  = \frac{{\theta _1  + 2u+ 2\theta \sqrt {\theta _2
} v (1 + u) + \theta_1u^2+ \theta ^2 \theta _1 \theta _2 v^2
}}{\theta^2\theta_1\theta_2+2\theta\sqrt\theta _2(u +v)+\theta _1
(u^2  +v^2)+ 2u v}\,.
\end{eqnarray}

To describe phases (Gibbs measures) of a given Hamiltonian on a
Cayley tree, one has a correspondence between Gibbs measures and a
collection of vectors $\{h_x\,, x \in V \}$, which satisfy a
non-linear equation. The recursive equations~(\ref{eq4})
considered in this paper describe a vector function
$\{(u_n\,,v_n), n \in \mathbb{N}\}$ which is a particular case of
the above mentioned function $h_x$ obtained as $h_x = u_n$ if $x
\in W_n$ i.e., depends only on the number of the generation set to
which belongs $x$ but not on $x$ itself (see for
example,~\cite{GTA, GR}).

The solutions of the system~(\ref{eq4}) describe the
translation-invariant Gibbs measures~\cite{P}. The number of the
solutions of the equations~(\ref{eq4}) depends on the parameters
$\beta={1}/{kT}$, $\theta$, $\theta _1$\,, $\theta_2$\,, $\theta
_3$ and $h$. The phase transition usually occurs for low
temperature. It is possible to find an exact value of $T^{\ast}$
such that a phase transition occurs for all $T<T^{\ast}$ where
$T^{\ast}$ is called a critical value of temperature.

An attractive fixed point of a function $f$ is a fixed point $u_0$
of $f$ such that for any value of $u$ in a domain that is close
enough to $u_0$\,, the iterated function sequence $u,\ f(u),\
f(f(u)),\ f(f(f(u))), \dots$ converges to $u_0$\,. An attractive
fixed point is said to be a stable fixed point if it is also
Lyapunov stable (see~\cite{Devaney} for details).

\section{Translation-invariant Gibbs measures}
As stated in~\cite{P}, for a given potential $V$ we have exactly
one Gibbs state with potential $V$ if the graph is finite. For
some potentials there may be more than one Gibbs state if the
graph is not finite. When we have more than one Gibbs state, then
we say that there exists a phase transition for a given potential
$V$. In other words, this occurrence of non-uniqueness of a Gibbs
measure can be interpreted as a phase transition.

In this section we determine whether any phase transition occurs
or not by solving the system~\eqref{eq4}. If there is more than
one positive solution for the equations~\eqref{eq4}, then there is
one Gibbs measure~\cite{MR} corresponding to each of these
solutions. If there is more than one Gibbs measure, then it is
said that a phase transition occurs for this model~\eqref{hm}
~\cite{GR,GTA}.

\subsection{First case}
Assume that $u=v$, then some solutions of the system~(\ref{eq4})
can be found from equation
\begin{equation}\label{eq8}
\theta_{3}u=\frac{\left(\tilde{\theta}^{2}\theta_{1}+2\tilde{\theta}+\theta_{1}\right)u^{2}
+2\left(\tilde{\theta}
+1\right)u+\theta_{1}}{2\left(\theta_{1}+1\right)u^{2}+4\tilde{\theta}
u+\tilde{\theta}^{2}\theta_{1}}:=f(u),
\end{equation}
where $\tilde{\theta}=\theta\sqrt {\theta _2 }$\,.

In this case, in order to obtain the solutions of the
system~\eqref{eq4}, we should analyze the equation in~\eqref{eq8}.
\begin{enumerate}
  \item When we take the first derivative of $f(u)$, we see
  that if $\tilde{\theta}>1$, that is $\theta _2
> {1}/{{\theta^2}}$, then $f(u)$ is increasing and if
$\tilde{\theta}<1$, $f(u)$ is
decreasing for $u>0$;
  \item When we take the second derivative of $f(u)$, we release that if
$\theta_{1}>\theta^{\ast}_{1}$\,, where
\[\theta^{\ast}_{1}=\frac{\tilde{\theta}^{2}+\tilde{\theta}+1+
\sqrt{9\tilde{\theta}^{4}+26\tilde{\theta}^{3}+35\tilde{\theta}^{2}+50\tilde{\theta}+33
}}{\left(\tilde{\theta}^{2}+2\right)\left(\tilde{\theta}+1\right)}\,,
\]
then there is an inflection point $u^{\ast}>0$ such that
$f''(u)>0$ for $0<u<u^{\ast}$ and $f''(u)<0$ for $u>u^{\ast}$
(see~\cite[Proposition 10.7]{P}  for details).
\end{enumerate}
\begin{theorem}\label{thm1}\cite[Theorem 5.1]{apostol} Suppose that
$f''$ exists on an interval $I$.

(i) If $f''(u)> 0$ on $I$, then the graph of $f$ is concave upward on $I$.

(ii) If $f''(u)<0$ on $I$, then the graph of $f$ is concave downward on $I$.
\end{theorem}
Thus we can generalize the Lemma proved in~\cite{GAT} to the
model~\eqref{hm} as follows.
\begin{lemma}
Assume that the quartic polynomial $g$ in~\eqref{eq10} has two
positive real roots. Let $\tilde{\theta}>1$, \linebreak $\theta _2
> {1}/{{\theta^2}}$, that is the equation $f$ in~\eqref{eq8} has a
positive inflection point. Then, there exist
$\eta_{1}(\theta,\theta_{1}\,,\theta_{2})$ and
$\eta_{2}(\theta,\theta_{1}\,,\theta_{2})$  with
$0<\eta_{1}(\theta,\theta_{1}\,,\theta_{2})<\eta_{2}(\theta,\theta_{1}\,,\theta_{2})$
such that the equation~\eqref{eq8} has three positive roots
$u_{\ast}^{(1)}<u_{\ast}^{(2)}<u_{\ast}^{(3)}$ if
$\eta_{1}(\theta,\theta_{1}\,,\theta_{2})<\theta_{3}<\eta_{2}(\theta,\theta_{1}\,,\theta_{2})$;
has two solutions if either
$\theta_{3}=\eta_{1}(\theta,\theta_{1}\,,\theta_{2})$ or
$\theta_{3}=\eta_{2}(\theta,\theta_{1}\,,\theta_{2})$ and has a
single solution in other cases. In fact
\[
\eta_{i}(\theta,\theta_{1}\,,\theta_{2})=\frac{1}{u^{\ast}_{i}}f(u^{\ast}_{i}),
\]
where $u^{\ast}_{1}\,, u^{\ast}_{2}$ are the solutions of equation
$uf'(u)=f(u).$
\end{lemma}
Proof follows from properties 1) and 2) of function $f$. As
mentioned in~\cite[Proposition 10.7]{P}, for the
equation~\eqref{eq8} there is more than one solution if and only
if there is more than one solution to $uf'(u)=f(u)$, which is
equivalent to the following equation
\begin{equation}\label{eq10}
g(u)=Au^{4}+Bu^{3}+Cu^{2}+Du+E=0,
\end{equation}
where $A=2\left(\theta _1+1\right)(\tilde{\theta}^{2}\theta
_1+2\tilde{\theta} + \theta _1)>0,$ $B=8(\tilde{\theta}
+1)\left(\theta _1+1\right)>0$, $D=8\tilde{\theta}
\theta_1>0$, $E=\tilde{\theta}^{2} \theta_1>0$ and
\[
C=-\left(\tilde{\theta}^{4}+\tilde{\theta}^{2}-6\right)\theta
^{2}_{1}-\left(2 \tilde{\theta} ^{3}-6\right)\theta
_1+8\tilde{\theta}^{2}.
\]

Let us determine where the graph of $g$ given in~\eqref{eq10} is
concave upward and concave downward. Here, we have
\[
g'(u)=4Au^{3}+3Bu^{2}+2Cu+D.
\]
We now have $g''(u)=12Au^{2}+6Bu+2C.$ Assume that
\begin{equation}\label{con}
3B^2-8AC>0. \
\end{equation}
From Theorem~\ref{thm1}, $g$ is concave downward on interval
\[
I=\left(\frac{-3B-\sqrt{9B^2-24AC}}{12A}\,,\,
\frac{-3B+\sqrt{9B^2-24AC}}{12A}\right)
\]
and  $g$ is concave upward on
interval $\mathbb{R}\setminus I$. Notice that at the point
\[
\left(\frac{-3B+\sqrt{9B^2-24AC}}{12A}\,,\,
g\left(\frac{-3B+\sqrt{9B^2-24AC}}{12A}\right)\right),
\]
the graph chances from concave downward to concave upward. Thus,
\[
\left(\frac{-3B+\sqrt{9B^2-24AC}}{12A}\, ,\,
g\left(\frac{-3B+\sqrt{9B^2-24AC}}{12A}\right)\right)
\]
is one of the inflection points of $g$.

\begin{wrapfigure}{o}{0.50\textwidth}
\centering\includegraphics[width=5cm]{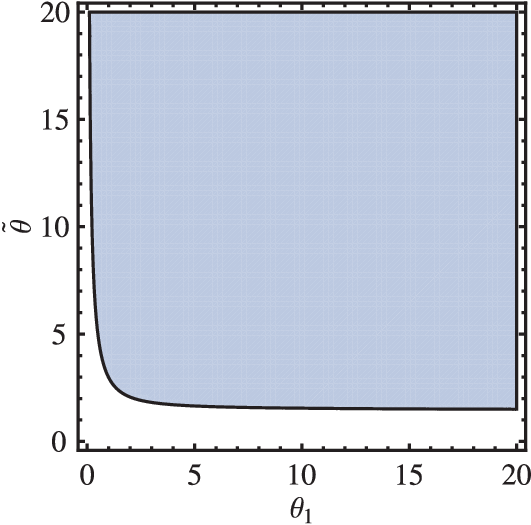}\
\caption{The region satisfying the inequality~\eqref{inequal}.}\label{region}%
\vspace{-5mm}
\end{wrapfigure}

It is clear that $g(0)=E>0$ and $\lim_{u\rightarrow\infty}g(u)=\infty$. Since $g$ is a continuous function, we can conclude that if  $g(u)<0$ on $J\subset \mathbb{R}^+$, then $g$ has two positive roots under condition \eqref{con}.

Due to Descartes' Rule of Signs~\cite{descard}, the maximum number
of positive real roots of the quartic polynomial $g$
in~\eqref{eq10} can be found by counting the number of sign
changes in the equation~\eqref{eq10}. Thus, the
equation~\eqref{eq10} has no positive roots if $C>0$, and the
polynomial~\eqref{eq10} has at most two positive roots
$u_{1}^{*}$\,, $u_{2}^{*}$ if $C<0$. In this case, if
$\tilde{\theta}>\sqrt{2}$ and $(2 \tilde{\theta}
^{3}-6)^{2}+32\tilde{\theta}^{2}(\tilde{\theta}^{4}+\tilde{\theta}^{2}-6)<0$
then the coefficient $C$ is negative.

Due to
\[
0=Au^4+Bu^3+Cu^2+Du+E\geqslant   Au^4+Cu^2+E\geqslant
2\sqrt{Au^4E}+Cu^2=(2\sqrt{AE}+C)u^2
\]
from this we get $2\sqrt{AE}+C\leqslant  0$ i.e. $C\leqslant
-2\sqrt{AE}$, this simple calculation shows that a necessary
condition to have positive roots of the equation \eqref{eq10} is
$C <-2\sqrt{AE}$. Thus, we have the following inequality:
\begin{equation}\label{inequal}
8 \tilde{\theta} + 8 \tilde{\theta}^2 + 6 \theta _1 - 2 \tilde{\theta}^3 \theta _1 + 6 \theta _1^2 - \tilde{\theta}^2 \theta _1^2 -\tilde{\theta}^4 \theta _1^2 +
  2 \sqrt{\tilde{\theta}^2 \theta _1 (1 + \theta _1)
  \left(2 \tilde{\theta} +\theta _1 +\tilde{\theta}^2 \theta _1\right)} < 0.
\end{equation}
Let us define the set $H=\left\{(\theta_1,\tilde{\theta}):C
<-2\sqrt{AE}\right\}$ (see figure~\ref{region}).

Then, we have the following
Corollary.
\begin{corollary}\label{lemma}
If $\theta _2
> {1}/{{\theta^2}}$, $3B^2-8AC>0$, $(\theta_1\,,\tilde{\theta})
\in H$ and $\eta_1(\theta,\theta_{1}\,,\theta
_2)<\theta_3<\eta_2(\theta,\theta_{1}\,,\theta _2)$ then the
function $f$  given in \eqref{eq8} has three fixed points. Thus, for the
model~\eqref{hm} there are three translation-invariant Gibbs
measures $\mu_1$\,, $\mu_2$\,, $\mu_3$ indicating a phase
transition.
\end{corollary}
It is easy to show that these three measures correspond to
boundary conditions $\bar{\sigma}(V\setminus V_{n})\equiv i,$ $
i=1, \ 2, \ 3.$

In figure~\ref{region} we show the region with three positive
roots. For parameters $\tilde{\theta}$, $\theta_1$ in phase
transition region the equation~\eqref{eq8} has three positive
roots, where two of them are stable and the third one is unstable.

\subsubsection{Illustrative example}
\begin{figure}[!h]
\centering
\includegraphics[width=0.31\columnwidth, clip]{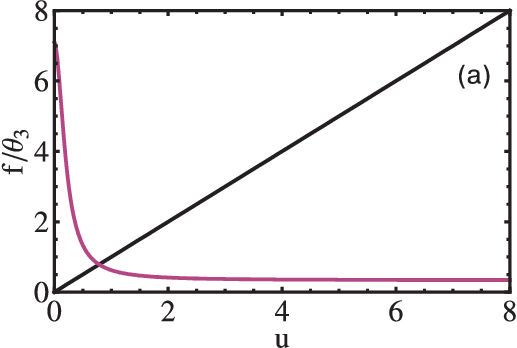}
\includegraphics[width=0.32\columnwidth, clip]{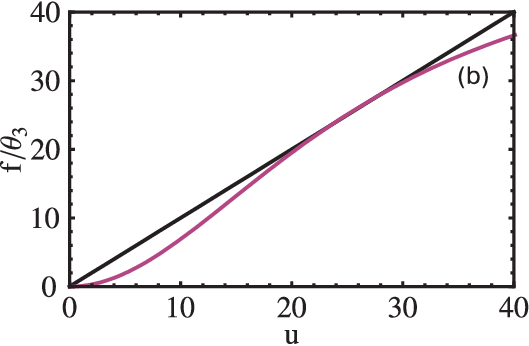}
\includegraphics[width=0.313\columnwidth, clip]{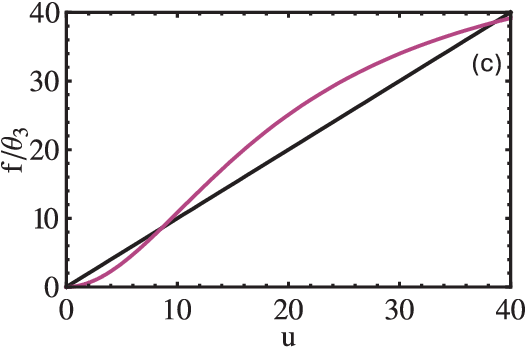}
\caption{Graphs of the function $f/\theta_{3}$ given
in~\eqref{eq8} (a) for $\theta=0.2$, $\theta_1=64$, $\theta
_{2}=2.2$ and $\theta_{3}=1.6$; (b) $\theta=11.6$, $\theta_1=28.1$,
$\theta _{2}=10.1$, $\theta_{3}=12.2$; (c) $\theta=10.3$,
$\theta_1=6.49$, $\theta _{2}=6.9$ and $\theta_{3}=6.4$,
respectively.}\label{exam1}
\end{figure}

Using the Mathematica software, we have manipulated the function
$f/\theta_{3}$ given in~\eqref{eq8} for some parameters $\theta$,
$\theta _{1}$\,, $\theta _{2}$ and $\theta_{3}$. In
figure~\ref{exam1}, for $\theta=0.2$, $\theta_1=64$, $\theta
_{2}=2.2$ and $\theta_{3}=1.6$ we have only one fixed point. Thus,
in this case the model~\eqref{hm} has no phase transition.

If we take as $\theta=11.6$, $\theta_1=28.1$, $\theta _{2}=10.1$,
$\theta_{3}=12.2$, then  we have two fixed points. As shown
figure~\ref{exam1}~(c), if we take as $\theta=10.3$, $\theta_1=6.49$,
$\theta _{2}=6.9$ and $\theta_{3}=6.4$, then  we have three fixed
points, where two of them are stable and the third one is
unstable.

Thus, in order to study the phase transition for the
model~\eqref{hm} we have clarified the role of $\theta$,
$\theta_1$\,, $\theta _{2}$ and $\theta_{3}$.


\subsection{Second case}
Now, let us find the other solutions of the system~\eqref{eq4}.
Assume that $u\neq v$. Subtracting the second equation
of~\eqref{eq4} from the first we have
\begin{equation}\label{eq7}
\theta_{3}(u-v)=\frac{2\left(\tilde{\theta}-1\right)(u-v)+\theta_{1}
\left(\tilde{\theta}^{2}-1\right)\left(u^{2}-v^{2}\right)}
{\tilde{\theta}^{2}\theta_{1}+2\tilde{\theta}(
u+v)+2uv+\theta_{1}\left(u^{2}+v^{2}\right)}\,.
\end{equation}
After canceling to $(u-v)$ two sides of the equation~\eqref{eq7}
we can obtain the following equation;
\begin{equation}\label{eq11}
t=\frac{\theta_{1}\theta_{2}s^{2}+\left[2\theta_{2}
\tilde{\theta}-\theta_{1}\left(\tilde{\theta}^{2}-1\right)\right]s
+\tilde{\theta}^{2}\theta_{1}
\theta_{2}-2\left(\tilde{\theta}-1\right)}
    {2\theta_{3}(\theta_{1}-1)}\,,
\end{equation}
where $u+v=s$  and  $uv=t$.

After dividing the first equation of~\eqref{eq4} to the second and
simplifying we have
\begin{equation}\label{eq12}
    t=\frac{\theta_{1}s^{2}+2s+\theta_{1}}
    {\tilde{\theta}^{2}\theta_{1}+\theta_{1}-2\tilde{\theta}}\,.
\end{equation}
From~\eqref{eq11} and~\eqref{eq12} we obtain the following
equality;
\begin{eqnarray}\label{eq13}
\nonumber
\frac{\theta_{1}\left(\tilde{\theta}+1\right)s+2}{\theta_{3}}&=&
\frac{\theta_{1}
\left[\theta_{1}\left(\tilde{\theta}+1\right)-2\right]s^{2}+2
\left[\theta_{1}\left(\tilde{\theta}^{2}+\tilde{\theta}-2\right)
-2\left(\tilde{\theta}+1\right)\right]s}
{\theta_{1}\left(\tilde{\theta}^{2}+1\right)-2\tilde{\theta}}\\
&&{}+ \frac{\theta_{1}
\left[\theta_{1}\left(\tilde{\theta}^{3}+\tilde{\theta}^{2}+2\tilde{\theta}+2\right)
-2\left(\tilde{\theta}^{2}+\tilde{\theta}+1\right)\right]}
{\theta_{1}\left(\tilde{\theta}^{2}+1\right)-2\tilde{\theta}}\,.
\end{eqnarray}
Let us consider a system of equations
\begin{eqnarray}\label{eq14}
 \nonumber u+v&=&s,\\
    uv&=&t.
\end{eqnarray}
From \eqref{eq14}, we can write $v=s-u$ and  $u(s-u)=t$, thus we have
\begin{equation}\label{polin}
u^2-su+t=0.
\end{equation}
To have two different real roots of the second-order equation
\eqref{polin}, inequality $s^2-4t>0$ should be valid. Then the
equations~\eqref{eq14} have solutions, if $t<{s^{2}}/{4}$.

Thus, from~\eqref{eq12} and from the inequality $t<{s^{2}}/{4}$ we
can get the following equation
\begin{equation}\label{eq15}
\left[\theta_{1}\left(\tilde{\theta}^{2}-3\right)-2\tilde{\theta}\right]s^{2}
-8s-4\theta_{1}>0.
\end{equation}
Then, there is
\[
s^{*}=\frac{4+\sqrt{16+4\theta_{1}\left[\theta_{1}\left(\tilde{\theta}^{2}-3\right)
-2\tilde{\theta}\right]}}
{\theta_{1}\tilde{\theta}^{2}-3\theta_{1}-2\tilde{\theta}}>0,
\]
such that for all $s>s^{*}$ the inequality $t<{s^{2}}/{4}$ is
valid.

One can show that the inequality~\eqref{eq15} is satisfied, if
\begin{equation}\label{eq15a}
\left(\tilde{\theta}^{2}-3\right)\theta_{1}^2-2\tilde{\theta}\theta_{1}+4<0
\end{equation}
and $\theta_{1}(\tilde{\theta}^{2}-3)-2\tilde{\theta}>0$.
Thus, for $\tilde{\theta}>\sqrt{3}$ and
\begin{equation}\label{eq15b}
0<\theta_{1}<\frac{\tilde{\theta}+\sqrt{12-3\tilde{\theta}^{2}}}
{\tilde{\theta}^{2}-3}
\end{equation}
the inequality~\eqref{eq15a} is valid.

In this case, the system of the equations~\eqref{eq14} has
solutions.

Assume that all coefficients of the equation~\eqref{eq13} are
positive. With the help of the elementary analysis we obtain the
following three cases:
\begin{enumerate}
\item If
\begin{equation}\label{eq16}
\theta _3 < \frac{{2\left[\theta _1 \left(\tilde{\theta} ^2  +
1\right) - 2\tilde{\theta} \right]}}{{\theta _1 \left[\theta _1
\left(\tilde{\theta} ^3  + \tilde{\theta }^2 + 2\tilde{\theta} +
2\right) - 2\left(\tilde{\theta} ^2  + \tilde{\theta}  +
1\right)\right]}}\,,
\end{equation}
then the equation~\eqref{eq13} has only one positive root. \item
If
\begin{equation}\label{eq17}
\theta_3> \frac{{2\left[\theta _1 \left(\tilde{\theta }^2  +
1\right) - 2\tilde{\theta} \right]}}{{\theta _1 \left[\theta _1
\left(\tilde{\theta }^3  + \tilde{\theta} ^2 + 2\tilde{\theta}  +
2\right) - 2\left(\tilde{\theta} ^2  + \tilde{\theta}  +
1\right)\right]}}
\end{equation}
there is $k_{0}$ such that the line $y-{2}/{\theta_3}=k_{0}s$ is a
tangent for the parabola on the right-hand side of~\eqref{eq13}
and then for $\theta_3<{{\theta _1 (\tilde{\theta} +
1)}}/{{k_{0}}}$ the equation~\eqref{eq13} has two positive roots
and only one of them is greater than $s^{*}$.

\item The equation~\eqref{eq13} has two positive roots, when
\[
\theta > \sqrt \frac{3}{\theta _2}\,,
\qquad
0<\theta_{1}<\frac{\tilde{\theta} +\sqrt{12-3\tilde{\theta}^{2}}}
{\tilde{\theta}^{2}-3}\,,
\qquad
\theta_{1}>\frac{2\tilde{\theta}}{\tilde{\theta}^{2}-3}
\]
and
\begin{equation}\label{eq18}
\frac{{2\left[\theta _1 \left(\tilde{\theta} ^2  + 1\right) -
2\tilde{\theta} \right]}}{{\theta _1 \left[\theta _1
\left(\tilde{\theta} ^3 + \tilde{\theta} ^2 + 2\tilde{\theta } +
2\right) - 2\left(\tilde{\theta }^2 + \tilde{\theta } +
1\right)\right]}} < \theta_3 < \frac{{\theta _1
\left(\tilde{\theta } + 1\right)}}{{k_0 }}\,,
\end{equation}
where $k_{0}$ is defined as above.
\end{enumerate}

\subsubsection{Illustrative example}
Let us give an illustrative example. With the help of the
Mathematica software, we have manipulated the line on the
left-hand side and the parabola on the right-hand side of the
equation~\eqref{eq13} for some parameters $\theta$, $\theta
_{1}$\,, $\theta _{2}$ and $\theta_{3}$. We have found the points
of intersection of the parabola and the line
given in~\eqref{eq13}.

\begin{figure}[h]
\centering
\includegraphics[width=0.31\columnwidth, clip]{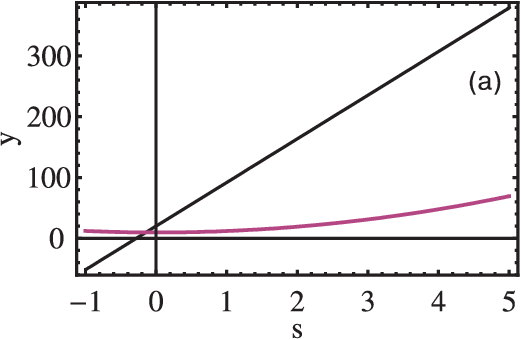}
\includegraphics[width=0.32\columnwidth, clip]{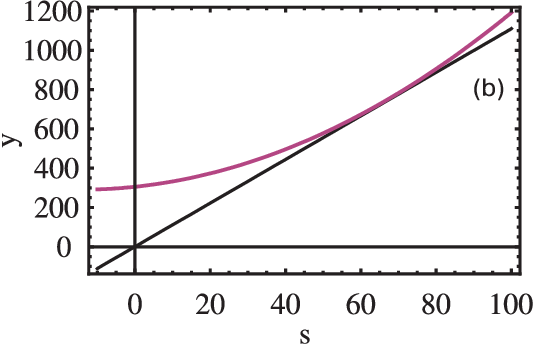}
\includegraphics[width=0.32\columnwidth, clip]{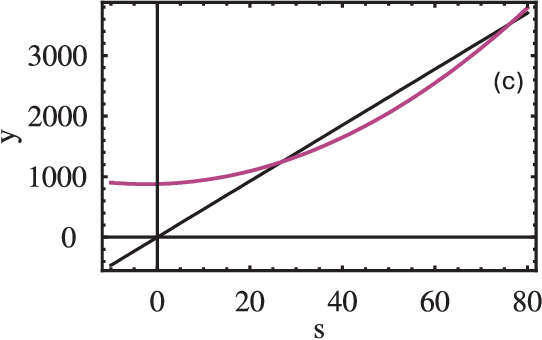}
\caption{Graphs of the parabola and the line given in~\eqref{eq13}
a) for $\theta=1.1$, $\theta_1=2.9$, $\theta _{2}=1.8$ and
$\theta_{3}=1.1$; b) $\theta=27.5$, $\theta_1=4.5$, $\theta
_{2}=5.9$ and $\theta_{3}=27.5$; c) $\theta=6.42$, $\theta_1=19$,
$\theta _{2}=50$ and $\theta_{3}=19$, respectively.}\label{exam2}
\end{figure}

In figure~\ref{exam2}~(a), the inequality~\eqref{eq16} is satisfied,
thus we have only one positive root. In figure~\ref{exam2}~(b), the line
$y-{2}/{\theta_3}={\theta_{1}(\tilde{\theta}+1)s}/{\theta_{3}}$ is
a tangent for the parabola~\eqref{eq13}. In figure~\ref{exam2}~(c), the
inequality~\eqref{eq18} is satisfied, we have two positive
solutions. Thus, in order to obtain the solutions of the
system~\eqref{eq14}, we have clarified the role of the parameters
$\theta$, $\theta_1$\,, $\theta _{2}$ and $\theta_{3}$\,.

Recalling the results for the equations~\eqref{eq4}, we get the following result:\\
if $(\theta, \ \theta_1\,, \ \theta _2)\in M$ and $K<\theta _3<L$,
where
\begin{eqnarray}\label{eq20}
\nonumber M&=&\left\{(\theta, \theta _1\,,\theta _2):\theta >
\sqrt \frac{3}{\theta _2}\,,
\
\theta_{1}>\frac{2\tilde{\theta}}
{\left(\tilde{\theta}^{2}-3\right)}\,,
\
0<\theta_{1}<\frac{\tilde{\theta}+\sqrt{12-3\tilde{\theta}^{2}}}
{\tilde{\theta}^{2}-3}
\right\}\cap H,\\
K&=& \max \left(\eta _1\,, \ \frac{{2\left[\theta _1
\left(\tilde{\theta} ^2 + 1\right) - 2\tilde{\theta}
\right]}}{{\theta _1 \left[\theta _1 \left(\tilde{\theta}^3 +
\tilde{\theta} ^2 + 2\tilde{\theta } + 2\right) -
2\left(\tilde{\theta }^2 +
\tilde{\theta } + 1\right)\right]}}\right),\nonumber\\
L&=& \min \left(\eta _2\,, \ \frac{{\theta _1
\left(\tilde{\theta}+1\right)}}{{k_{0}}}\right)
\end{eqnarray}
then the system~\eqref{eq4} has five positive solutions, three of
which are obtained from the first case and the other two solutions
arise in the second case. Similar to~\cite{GAT} and~\cite{G}, a
more detailed analysis shows that only three solutions are stable.

Combining the first and the second cases, we have proved the following theorem.
\begin{theorem} Assume that the conditions~\eqref{eq20} are satisfied,
then for the model~\eqref{hm} there are five translation-invariant
Gibbs measures indicating a phase transition.
\end{theorem}

\begin{rmk} In a similar way we can study the solutions of a new system of
nonlinear equations by using Kroneker's symbol $\delta$ that has
the form
$$
\delta _{\sigma (x)\sigma (y)\sigma (z)} =\left\{
{{\begin{array}{*{20}c}
 1, \hfill & \text{if}\ \ {\sigma (x)=\sigma (y)=\sigma (z)} \hfill, \\
 0, \hfill & \text{otherwise}. \hfill \hfill \\
\end{array}}} \right.
$$
\end{rmk}

\section{Conclusion}

In this paper we have exactly solved the Potts model on a Cayley
tree, the Hamiltonian of which contains three competing
interactions, the nearest neighbors, the second neighbors, triples
of neighbors  and the external field. Namely, we have calculated
the critical curve such that there is a phase transition above it,
and a single Gibbs state is found elsewhere. We have clarified the
role of the coupling constants $J$, $J_{1}$\,, $J_{2}$ and
external field $h$ to study the existence problem of a phase
transition. We have considered transformation $\textbf{F} =
(F_1\,; F_2) : \mathbb{R}^{2} \rightarrow \mathbb{R}^{2}$ and
proved that its fixed points describe the translation-invariant
Gibbs measures of our model. We have seen that for some values
$J$, $J_{1}$\,, $J_{2}$ and $h$, the phase transition occurs.

In~\cite{GR}, an Ising model with four competing interactions
(external field, nearest neighbor, second neighbors and triples of
neighbors) on the Cayley tree of order two has been solved.  The
authors~\cite{GR} have also constructed numerous non-periodic
extreme Gibbs measures. The investigation of the periodic Gibbs
measures for the Potts model~\eqref{hm} is planned to be the
subject of forthcoming publications.

\section*{Acknowledgement}
The authors would like to thank Prof. Dr.~N.N.~Ganikhodjaev and
Prof. Dr.~Utkir Rozikov for valuable support and suggestions. We
are also  grateful to anonymous referees whose comments greatly
improved the paper.

\ukrainianpart

\title%
{До фазових переходів моделі Поттса з трьома конкуруючими взаємодіями на дереві Келі
}

\author[Г. Акін, С. Темір]
{Г.~Акін\refaddr{label1}, С.~Темір\,\refaddr{label2}}

\addresses{
\addr{label1} Університет Зірве, 27260  Газіантеп, Туреччина
\addr{label2} Університет Гарран,  63200 Санліурфа, Туреччина}

\makeukrtitle

\begin{abstract}
У цій статті ми вивчаємо проблему фазового переходу  для моделі Поттса з трьома конкуруючими взаємодіями, найближчих сусідів,
наступних  близьких сусідів  і наступних за наступними близькими сусідами, з ненульовим зовнішнім полем на дереві Келі другого порядку.
Ми доводимо, що для деяких значень параметрів моделі є фазовий перехід. Ми приводимо проблему опису граничних мір Гіббса до проблеми  розв'язку
системи нелінійних функціональних рівнянь. Ми розширюємо результати, отримані Ганіходяєвим і Розіковим
[Math.
Phys. Anal. Geom., 2009, \textbf{12}, No~2, 141--156]
для фазового переходу в моделі Ізінга, на випадок
моделі Поттса.

\keywords фазовий перехід, модель Поттса, конкуруючі взаємодії, міра Гіббса
\end{abstract}

\label{last@page}

\end{document}